 \RequirePackage{fix-cm}
\documentclass[smallextended, envcountsame]{svjour3}
\smartqed
\def\makeheadbox{\relax}
\journalname{}

\usepackage[unicode = false,
    pdftoolbar = true,
    colorlinks = true,
    linkcolor = blue,
    citecolor = blue,
    filecolor = black,
    urlcolor = blue,
    breaklinks = true]{hyperref}

\usepackage[utf8]{inputenc}
\usepackage{listings}  
\usepackage{enumitem}  
\usepackage{amsmath}   
\usepackage{amssymb}   
\usepackage{xcolor}    
\usepackage{graphicx}  
\usepackage{verbatim}
\usepackage{geometry}
\usepackage{mathtools}
\usepackage{float}
\usepackage{bm}
\usepackage[ruled]{algorithm2e}
\usepackage{subcaption}
\usepackage{placeins}
\usepackage{listings}
\usepackage[justification=centering]{caption}
\usepackage{multirow}
\usepackage{etoolbox}
\usepackage{orcidlink}
\usepackage{bbding}
\AtEndEnvironment{proof}{\leavevmode\qed}

\usepackage{soul}
\definecolor{gab}{HTML}{50c878}
\definecolor{pat}{HTML}{bf68ff}
\definecolor{sham}{HTML}{7faaff}

\allowdisplaybreaks[4]

\newtheorem{theo}{Theorem}

\DeclareMathAlphabet{\mymathbb}{U}{BOONDOX-ds}{m}{n}

\makeatletter
\newcommand*{\rom}[1]{\expandafter\@slowromancap\romannumeral #1@}
\makeatother

\SetKwIF{If}{ElseIf}{Else}{if}{}{else if}{else}{end if}%
\SetKwRepeat{Do}{do}{while}
\allowdisplaybreaks[4]
\SetArgSty{textnormal}

\title{Non-realizability of 4091 topological prismatoids}
\author{Yiwen Chen\,\orcidlink{0000-0002-8720-932X}\,\Envelope\and Amy Wiebe\,\orcidlink{0000-0002-0804-6252}}
\institute{Department of Mathematics, University of British Columbia, Kelowna, British Columbia, V1V 1V7, Canada.\\
This research is partially funded by the Natural Sciences and Engineering Research Council of Canada (cette recherche est partiellement financ\'ee par le Conseil de recherches en sciences naturelles et en g\'enie du Canada), Discover Grant \#2024-04643.\\
\email{yiwchen@student.ubc.ca, amy.wiebe@ubc.ca}}
\date{\today}

\begin{document}

\maketitle

\begin{abstract}
    Motivated by the search for small counterexamples to the Hirsch conjecture, Criado and Santos constructed a collection of topological prismatoids whose associated simplicial spheres violate the Hirsch bound.  If any of the sufficiently small examples admits a convex realization, then it could lead to a counterexample in lower dimensions than those previously constructed. However, the realizability of these prismatoids was left as an open question. Subsequent work by Pfeifle and by Gouveia, Macchia, and Wiebe established non-realizability for subsets of the collection. In this paper, we use linear programming to show that all 4091 prismatoids from Criado and Santos\footnote{Although Criado and Santos mentioned 4093 prismatoids in their original work \cite{criado2022topological}, their published data contains only 4091.  In personal communication, they have stated that they indeed only have 4091  prismatoids in the collection.} are non-realizable. By adapting the computational framework of Gouveia, Macchia, and Wiebe, we obtain exact non-realizability certificates for every prismatoid in the collection.  Remarkably, these certificates require only products of at most two polynomials with each of degree at most three. Our results rule out all 4091 prismatoids as potential sources of convex counterexamples to the Hirsch conjecture, closing the realizability question for this collection.
\end{abstract}
\medskip

\noindent {\bf Keywords:} Non-realizability certificates; Hirsch conjecture; Linear programming
\medskip

\noindent {\bf MSC codes:} 52B11; 52B05; 52B55; 90C05

\section{Introduction}
\label{sec:introduction}

    The Hirsch conjecture asserted that the graph of a convex $d$-polytope with $m$ facets has diameter at most $m-d$ \cite{dantzig1965linear,santos2012counterexample}. Here, the graph consists of the vertices and edges of the polytope, and its diameter is the largest number of edges in a shortest path between two vertices. The conjecture was motivated by linear programming: the simplex method moves between adjacent vertices of a feasible polyhedron, so distances in this vertex-edge graph give lower bounds on the number of steps needed to reach an optimal vertex \cite{dantzig1965linear}. We say that a polytope is \emph{non-Hirsch} if its graph diameter exceeds $m-d$.

    Santos~\cite{santos2012counterexample} disproved the conjecture by constructing a non-Hirsch polytope in dimension 43. Matschke, Santos, and Weibel~\cite{matschke2015width} subsequently reduced the dimension to 20. Both constructions begin with a \emph{prismatoid}, a convex polytope whose vertices lie in two parallel facets, called its \emph{bases}. The \emph{width} of a prismatoid is the distance between its two bases, measured by the minimum number of steps between adjacent facets.  When the width exceeds the dimension of the prismatoid, the prismatoid can be used to construct a higher-dimensional polytope whose graph diameter violates the Hirsch bound. A starting prismatoid with fewer vertices can lead to a counterexample in lower dimensions.

    Recent work has also studied how large the width of a prismatoid and the graph diameter of a polytope can be. Davies et al.~\cite{davies2025drums} constructed a family of 5-dimensional prismatoids whose width grows linearly in the number of vertices. More recently, OpenAI~\cite{openai2026superpolynomial} reported a disproof of the polynomial Hirsch conjecture by constructing a family of polytopes whose graph diameter grows faster than any fixed power of the number of facets. In contrast, the problem that motivates this paper is to find prismatoids with few vertices that could lead to counterexamples to the original Hirsch conjecture in lower dimensions.

    Criado and Santos~\cite{criado2022topological} searched for smaller examples in a broader class of objects called \emph{topological prismatoids}. These are combinatorial complexes that model the boundary of a prismatoid with its two bases removed, without requiring a convex realization. Their search produced 4091 topological prismatoids of dimension 4 with the required distance property, ranging from 14 to 28 vertices. If one of the smallest examples could be realized as a convex 5-dimensional prismatoid, it would yield a non-Hirsch polytope in dimension 9. 

    Pfeifle~\cite{pfeifle2024positive} and Gouveia, Macchia, and Wiebe~\cite{gouveia2023general} subsequently established non-realizability for subsets of these prismatoids. In this paper, we apply the method of~\cite{gouveia2023general}, using a modified version of their code, to all 4091 prismatoids from~\cite{criado2022topological}.

    \begin{theo}\label{thm:nonrealizability}
        All 4091 topological prismatoids from the search of Criado and Santos~\cite{criado2022topological} are non-realizable. More precisely, after adjoining the two bases to each input facet description, none of the resulting vertex-facet incidence structures is realizable by a convex 5-polytope.
    \end{theo}

    Thus none of these candidates can be realized as a convex prismatoid from which to construct a counterexample to the Hirsch conjecture. 
    
    The remainder of this paper is organized as follows. Section~\ref{sec:background} introduces the collection of 4091 topological prismatoids constructed in~\cite{criado2022topological}. Section~\ref{sec:certificates} describes the non-realizability certificates introduced in~\cite{gouveia2023general}.  Section~\ref{sec:results} reports our computational results.  Section~\ref{sec:conclusion} summarizes this paper and discusses future directions.

\section{Background and the prismatoid collection}
\label{sec:background}

\subsection{Geometric prismatoids and the strong $d$-step theorem}

    The constructions of the prismatoids discussed in this section rely on distances between facets rather than vertices. This motivates considering the \emph{facet adjacency graph}, or dual graph, of a $d$-polytope, which has one node for each facet, with two nodes adjacent when the corresponding facets intersect in a face of dimension $d-2$. Polarity exchanges vertices and facets, identifying the dual graph with the vertex-edge graph of the polar polytope. Consequently, the dual form of the Hirsch bound for a $d$-polytope with $n$ vertices states that its facet adjacency graph has diameter at most $n-d$.

    A geometric $d$-prismatoid is a $d$-polytope with two parallel facets containing all its vertices.  Recall that a convex polytope is called \emph{simplicial} if all its facets are simplices.  For a geometric prismatoid, we require only the facets other than the two bases to be simplices for it to be called simplicial.  The \emph{width} $\ell$ of a geometric $d$-prismatoid is the length of a shortest path between the two bases in the facet adjacency graph. If $\ell>d$, then we call the prismatoid \emph{non-$d$-step}.  Santos' strong $d$-step theorem~\cite[Theorem~2.6]{santos2012counterexample} converts this width condition into a violation of the Hirsch bound. In particular, given a simplicial geometric $d$-prismatoid with $n>2d$ vertices and width $\ell$, the theorem constructs a geometric prismatoid of dimension $D=n-d$ with $N=2(n-d)$ vertices and width at least $\ell+n-2d$. If $\ell>d$, then $\ell+n-2d>n-d=N-D$. Hence, the constructed prismatoid violates the dual Hirsch bound, and its polar is a non-Hirsch $D$-polytope with $N$ facets.

    For fixed $d$, reducing the number $n$ of vertices of the starting prismatoid reduces the dimension $D=n-d$ of the resulting counterexample. Santos' original construction used a 5-prismatoid with 48 vertices, giving dimension 43~\cite{santos2012counterexample}. The 25-vertex example of Matschke, Santos, and Weibel gives dimension 20~\cite{matschke2015width}. These constructions motivate the search for smaller non-$d$-step prismatoids.

\subsection{Topological prismatoids and non-Hirsch spheres}
    We now extend the geometric construction discussed in the previous subsection to topological prismatoids and non-Hirsch spheres.  We first introduce simplicial complexes.  Recall that, for a simplicial polytope, every face is a simplex, so every subset of the vertices of a face also forms a face.  A \emph{simplicial complex} generalizes this combinatorial structure: it is a collection of subsets of a vertex set that is closed under taking subsets. A \emph{simplicial $(d-1)$-sphere} is a simplicial complex that triangulates the $(d-1)$-dimensional sphere $S^{d-1}$. A simplicial $(d-1)$-sphere is \emph{realizable} if it is combinatorially isomorphic to the boundary complex of a convex $d$-polytope.

    Notice that a simplicial geometric $d$-prismatoid gives rise to a simplicial complex by removing the relative interiors of its two bases from its boundary. The resulting complex is homeomorphic to $S^{d-2}\times[0,1]$, a cylinder whose two boundary components correspond to the boundaries of the bases~\cite{criado2022topological}. A \emph{topological $(d-1)$-prismatoid} generalizes this structure: it is a simplicial complex of dimension $d-1$ that is homeomorphic to $S^{d-2}\times[0,1]$ and whose two boundary components are \emph{induced subcomplexes}~\cite{criado2022topological}. Here, induced means that any face whose vertices all lie in one boundary component also belongs to that component.  We say that such a complex admits a \emph{convex realization} if it is combinatorially isomorphic to the complex obtained from the boundary of a convex $d$-polytope by removing the relative interiors of two disjoint facets containing all its vertices.

    We can naturally extend the notion of width from geometric $d$-prismatoids to topological $(d-1)$-prismatoids. As before, we consider the facet adjacency graph, where each node represents a maximal face of the simplicial complex, and two nodes are adjacent if the corresponding faces share a $(d-2)$-dimensional face. The \emph{width} of a topological $(d-1)$-prismatoid is defined as two plus the minimum distance in this graph between a facet adjacent to one boundary component and a facet adjacent to the other. The additional two steps account for entering from one base and exiting to the other when the bases are restored. Thus, this definition agrees with geometric width whenever a convex realization exists. A topological $(d-1)$-prismatoid is called \emph{non-$d$-step} if its width exceeds $d$.

    The significance of non-$d$-step topological prismatoids follows from the topological strong $d$-step theorem~\cite[Theorem~3.5]{criado2022topological}, which extends Santos' strong $d$-step theorem~\cite[Theorem~2.6]{santos2012counterexample} from geometric to topological prismatoids. Specifically, if a non-$d$-step topological $(d-1)$-prismatoid has $n>2d$ vertices and both boundary components are realizable as boundaries of convex polytopes, the theorem yields a simplicial $(D-1)$-sphere with $N=2D$ vertices and facet adjacency diameter greater than $N-D$, where $D=n-d$. A simplicial $(D-1)$-sphere with $N$ vertices is called \emph{non-Hirsch} if its facet adjacency diameter exceeds $N-D$. If the resulting sphere were realizable as the boundary of a convex $D$-polytope, its polar would violate the Hirsch bound.

\subsection{The collection and existing results}

    Criado and Santos~\cite{criado2022topological} used simulated annealing to search for small non-$d$-step topological 4-prismatoids. Their search employed local modifications of simplicial complexes that preserve their topology, while maintaining the non-$d$-step condition and favoring reductions in the number of vertices. However, these modifications need not preserve convex realizability. Thus, even when starting from a geometric prismatoid, the resulting topological prismatoids need not be realizable.
    
    Their search started from the 28-vertex geometric prismatoid described in~\cite{matschke2015width}. They reported 4091 completed experiments, producing topological prismatoids with between 14 and 28 vertices~\cite[Section~5]{criado2022topological}. Notice that this collection is a computational sample rather than an exhaustive enumeration of all topological prismatoids in this vertex range.

    The realizability of these topological prismatoids was left open in~\cite{criado2022topological}. Pfeifle~\cite{pfeifle2024positive} subsequently proved that the four combinatorially distinct 14-vertex examples in the collection are non-realizable using positive Pl\"ucker tree certificates. These certificates combine Grassmann-Pl\"ucker relations \cite[Section~1.2]{bjorner1999oriented} to obtain expressions that vanish algebraically but would be strictly positive under a convex realization.

    Gouveia, Macchia, and Wiebe~\cite{gouveia2023general} developed a method based on slack matrices and linear programming. A slack matrix records the slack of each vertex in each facet inequality (see Section~\ref{sec:certificates} for details). Using this approach, they obtained non-realizability certificates for all four examples with 14 vertices, all 36 examples with 15 vertices, and selected larger examples, including the unique 28-vertex prismatoid in the collection. They raised the questions of whether all prismatoids in the collection are non-realizable and whether a common obstruction could explain their non-realizability. Pfeifle's subsequent work~\cite{pfeifle23proving} provided further certificates based on related determinant identities, but the realizability of the full collection remained unresolved. This paper resolves the non-realizability question for the entire collection by establishing certificates for all 4091 prismatoids, as presented in Theorem~\ref{thm:nonrealizability}.

\section{Non-realizability certificates}
\label{sec:certificates}

    We use the method of Gouveia, Macchia, and Wiebe~\cite{gouveia2023general} to construct polynomials that must be positive in any realization and then search for a nonnegative combination that is identically zero.  In this section, we recall the notation needed to describe our certificates and the bounds used in the search.

\subsection{Slack matrices and positive generators}

    Let $P$ denote the vertex-facet incidence structure resulting from adding the two bases to the facet list of a topological 4-prismatoid. The \emph{symbolic slack matrix} $S_P(z)$ has one row for each vertex $i$ and one column for each facet $F$, where $F$ is identified with its set of vertices:
    \begin{equation*}
        (S_P(z))_{i,F}=
        \begin{cases}
            0,&i\in F,\\
            z_{i,F},&i\notin F,
        \end{cases}
    \end{equation*}
    where the $z_{i,F}$ are distinct variables. In a convex realization with vertices $v_i\in\mathbb R^5$ and facet inequalities $a_F^\top v\le b_F$, the entries of the slack matrix are $b_F-a_F^\top v_i$. They vanish when $i\in F$ and are strictly positive otherwise.  Notice that a slack matrix of a convex 5-polytope has rank six~\cite{gouveia2013nonnegative}, whereas an arbitrary substitution into $S_P(z)$ need not satisfy this rank condition.

    Following~\cite[Section~4]{gouveia2023general}, we choose six columns of the symbolic slack matrix, including those corresponding to the two bases, that are necessarily linearly independent in any convex realization. These columns form a \emph{reduced slack matrix}.  To reduce the number of variables, we use positive row and column scalings to set certain nonzero entries to one. Specifically, we choose entries corresponding to a spanning tree of the bipartite graph associated with the nonzero pattern of the reduced slack matrix, as described in~\cite[Section~4]{gouveia2023general}. Let $x$ denote the remaining variables, and let $u_i(x)$ be row $i$ of the scaled reduced slack matrix.

    In any convex realization, the rows $u_i(x)$ represent homogeneous vertex coordinates.  For each facet $F$, we choose five vertices $(b_1,\ldots,b_5)$ whose affine hull would be the supporting hyperplane of $F$ in any realization. We call this choice a \emph{facet basis}. If $F$ is simplicial, its five vertices form a facet basis. For a nonsimplicial facet, we choose a suitable set using the face incidences, as described in~\cite[Section~2]{gouveia2023general}. We define 
    \begin{equation*}
        p_{i,F}(x)=\varepsilon_F
        \det\begin{pmatrix}
            u_{b_1}(x)\\
            \vdots\\
            u_{b_5}(x)\\
            u_i(x)
        \end{pmatrix}~~~\text{where}~~~\varepsilon_F\in\{-1,1\}.
    \end{equation*}
   With the facet basis fixed, the determinant is linear in its last row and vanishes on the supporting hyperplane of $F$. Its sign therefore distinguishes the two sides of that hyperplane.  All vertices outside a facet of a convex polytope lie strictly on the same side of that hyperplane. Consequently, the determinants for $i\notin F$ have the same nonzero sign. We choose $\varepsilon_F$ so that $p_{i,F}(x)>0$ for every $i\notin F$ in any realization.

    To determine $\varepsilon_F$, we follow the procedure in~\cite[Section~5.3.1]{gouveia2023general}. Since the variables in $x$ are positive, the sign of a determinant that is a nonzero monomial is determined by its coefficient. More generally, suppose a determinant before multiplication by $\varepsilon_F$ factors as $c\,q_1(x)\cdots q_r(x)$, where $c\ne0$ and each $q_k$ is already known to be positive in any realization. The determinant then has the sign of $c$. One such determinant fixes the sign for the entire facet, and the resulting positive expressions can help determine signs for other facets.

    We simplify the signed determinants $p_{i,F}(x)$ for $i\notin F$ by removing common monomial factors and dividing by the positive greatest common divisor of their integer coefficients. These operations preserve positivity because all variables in $x$ are positive. The simplified polynomials, together with the variables in $x$ and the constant $1$, form the \emph{generators} for our search. We retain generators of degree at most three and consider selected products of at most two generators.

\subsection{Certificates and exact verification}

    Let $g_1(x),\ldots,g_m(x)$ denote the products selected above. Each is strictly positive in any convex realization. We seek a polynomial identity of the form
    \begin{equation}\label{eq:certificate}
        \sum_{j=1}^{m}\lambda_j g_j(x)\equiv0,
        \qquad \lambda_j\ge0,
        \qquad \sum_{j=1}^{m}\lambda_j=1.
    \end{equation}
    Here, $\equiv0$ means that the polynomial vanishes identically, not just at a particular choice of $x$. Since the weights $\lambda_j$ are nonnegative and sum to one, the left-hand side would be strictly positive in any realization. Thus, an identity of this form proves non-realizability.
    
    Following~\cite[Section~5.1]{gouveia2023general}, we search for these identities by equating the coefficients of corresponding monomials in~\eqref{eq:certificate}. This gives a linear feasibility problem in the weights $\lambda_j$ that can be solved numerically. A candidate certificate can then be checked in exact rational arithmetic, so its validity does not depend on numerical tolerances.

\section{Computational results}
\label{sec:results}

    The original input data lists the facets other than the two bases. For each prismatoid, we recovered the bases from the vertex labels, relabeled the vertices consecutively, and applied the method described in Section~\ref{sec:certificates}. Using products of at most two generators, each of degree at most three, we obtained an exact non-realizability certificate for all 4091 prismatoids. Table~\ref{tab:counts} gives the number of prismatoids for each vertex count.

    \begin{table}[htbp]
    \caption{Distribution of the 4091 prismatoids by number of vertices.}
    \label{tab:counts}
    \centering
    \begin{tabular}{rr@{\qquad}rr@{\qquad}rr}
    \hline
    Vertices & Prismatoids & Vertices & Prismatoids & Vertices & Prismatoids\\
    \hline
    14 &   4 & 19 & 665 & 24 & 176\\
    15 &  36 & 20 & 695 & 25 &  76\\
    16 & 209 & 21 & 634 & 26 &  26\\
    17 & 372 & 22 & 435 & 27 &   7\\
    18 & 481 & 23 & 274 & 28 &   1\\
    \hline
    \end{tabular}
    \end{table}

    The number of nonzero summands in~\eqref{eq:certificate} ranges from three to eighteen. In particular, of the 4091 certificates, 1946 have three summands, and 1276 have four.

    \begin{proof}[of Theorem~\ref{thm:nonrealizability}]
        For each prismatoid, we verified in exact rational arithmetic that the computed coefficients satisfy~\eqref{eq:certificate}. If a convex realization existed, every $g_j(x)$ would be strictly positive by the construction in Section~\ref{sec:certificates}. Since the coefficients are nonnegative and sum to one, their weighted sum would also be strictly positive, contradicting~\eqref{eq:certificate}.
    \end{proof}

    The original input data, modified search code, and exact certificates are available in the repository~\cite{Data}.

\section{Conclusion}
\label{sec:conclusion}
    We have shown that none of the 4091 topological prismatoids found by Criado and Santos~\cite{criado2022topological} admits a convex realization, answering the question raised in~\cite{gouveia2023general}. Thus, none can serve as a convex starting point for the strong $d$-step construction. However, the collection is a finite computational sample, so our result does not rule out other realizable non-$d$-step prismatoids or establish a general lower bound on the size of non-Hirsch polytopes. 

    It remains to understand whether the certificates reflect common combinatorial obstructions to realizability. Identifying such obstructions could help explain why these examples are non-realizable and guide future searches for realizable non-$d$-step prismatoids.

\normalsize
\bibliographystyle{siam}
\bibliography{references}
\end{document}